\documentclass[graybox]{svmult}

\usepackage{mathptmx}       % selects Times Roman as basic font
\usepackage{helvet}         % selects Helvetica as sans-serif font
\usepackage{courier}        % selects Courier as typewriter font
\usepackage{type1cm}        % activate if the above 3 fonts are
\usepackage{makeidx}         % allows index generation
\usepackage{graphicx}        % standard LaTeX graphics tool
\usepackage{multicol}        % used for the two-column index
\usepackage[bottom]{footmisc}% places footnotes at page bottom
\usepackage{amssymb}
\usepackage{amsfonts}
\usepackage{amsmath}
\usepackage{graphicx}
\usepackage{mathtools}
\mathtoolsset{showonlyrefs}
\usepackage{tikz}
\usetikzlibrary{matrix}

\newcommand{\R}[1]{\mathbb{R}^{#1}}
\renewcommand{\S}[1]{\mathbb{S}^{#1}}

\newcommand{\cA}{\mathcal A}

\newcommand{\cH}{\mathcal H}

\newcommand{\cL}{\mathcal L}
\newcommand{\cM}{\mathcal M}

\newcommand{\cV}{\mathcal V}

\newcommand{\de}{\mathrm d}

\newcommand{\cl}[1]{\overline{#1}}

\newcommand{\eps}{\varepsilon}

\newcommand{\curl}{\operatorname{curl}}

\newcommand{\tr}{\operatorname{tr}}

\renewcommand{\leq}{\leqslant}
\newcommand{\wto}{\rightharpoonup}
\newcommand{\wsto}{\stackrel{*}{\rightharpoonup}}

\newcommand{\average}{{\mathchoice {\kern1ex\vcenter{\hrule
height.4pt width 8pt depth0pt}
\kern-11pt} {\kern1ex\vcenter{\hrule height.4pt width 4.3pt
depth0pt} \kern-7pt} {} {} }}

\newcommand{\res}{\mathop{\hbox{\vrule height 7pt width .5pt depth
0pt\vrule height .5pt width 6pt depth0pt}}\nolimits}

\makeindex             % used for the subject index
\begin{document}

\title*{Structured Deformations of Continua: Theory and Applications}
\author{Marco Morandotti}
\institute{Marco Morandotti \at Technische Universit\"at M\"unchen,  \email{marco.morandotti@ma.tum.de}}
\maketitle

\abstract*{The scope of this contribution is to present an overview of the theory of structured deformations of continua, together with some applications.
Structured deformations aim at being a unified theory in which elastic and plastic behaviours, as well as fractures and defects can be described in a single setting. 
Since its introduction in the scientific community of rational mechanicists \cite{DPO1993}, the theory has been put in the framework of variational calculus \cite{CF1997}, thus allowing for solution of problems via energy minimization.
Some background, three problems and a discussion on future directions are presented.}

\abstract{The scope of this contribution is to present an overview of the theory of structured deformations of continua, together with some applications.
Structured deformations aim at being a unified theory in which elastic and plastic behaviours, as well as fractures and defects can be described in a single setting. 
Since its introduction in the scientific community of rational mechanicists \cite{DPO1993}, the theory has been put in the framework of variational calculus \cite{CF1997}, thus allowing for solution of problems via energy minimization.
Some background, three problems and a discussion on future directions are presented.
}

\section{Introduction}
\label{sec:1}
Many and complex are the deformations of a body when it responds to external loading.
In the last century, and even before, several theories have been proposed to study phenomena such as elasticity, plasticity, fracture\ldots 
More recently, both the advent of modern instruments capable of resolving the core of the matter, and more powerful computers allowed scientists to bridge different length scales together and to formulate theories that range from the subatomic level to the macroscopic one.
In this way, phenomena that were described from a macroscopic viewpoint started being related to their microscopic counterpart: it was possible to model at the atomistic level and to explain phenomena that are experienced at the macroscopic scale.
An example of this is the understanding that the motion and pile-up of defects such as dislocations is the microscopic process which is responsible for plasticity.

The necessity of having a theoretical tool that permits inclusion of multiple scales is evident.
\emph{Structured deformations} \cite{DPO1993} (see also \cite{BMS2012}) respond to this need by providing a multiscale geometry that captures the contributions at the macroscopic level of both smooth geometrical changes and non-smooth geometrical changes at submacroscopic levels.
These non-smooth geometrical changes, which are called \emph{disarrangements}, encode the presence of cracks and defects in the continuum.
The geometrical effects at the macrolevel of submicroscopic disarrangements are captured by sequences of approximating, piecewise smooth deformations $f_n$ that converge to the macroscopic deformation field $g$ and whose gradients $\nabla f_n$ converge to a field $G$ that might differ from $\nabla g$.
The jumps in $f_n$ contribute to the interfacial part of an initial energy response $E(f_n)$, they diffuse throughout portions of the body, and they contribute in the limit to both the bulk and interfacial parts of a relaxed energy response $I(g,G)$.
%In the mathematical model proposed in \cite{DPO1993} and then set in a variational framework in \cite{CF1997} , the singularities of the deformation map, in principle entering the energy functional in the surface part, can diffuse to contribute to the bulk energy.
This is the main result contained in the approximation and relaxation theorems \cite[Theorem 5.8]{DPO1993} and \cite[Theorems 2.12 and 2.17]{CF1997}, and is the novelty of the theory.
The main feature of the relaxation theorems is that they provide a representation formula for the relaxed energy functional which can then be expressed as an integral.
The relaxed bulk and surface energy densities are obtained through the so-called cell formulas, which are expressed as minimum problems.

Structured deformations have been successfully applied in many contexts to model the deformation of bodies that also include plastic deformations and cracks \cite{CDPFO1999,DPT2001,DO2003,DO2015}.
The theory has been extended to more general contexts, especially by defining second-order structured deformations \cite{OP2000}, which permit the inclusion of bending effects in the energy functional.
Applications and developments of second-order structured deformations are provided, for example, in \cite{BMMO,O,P2004}.
Also relevant are the works \cite{BMMO-expl,OP2015,S}, which focus on interfacial energies, relevant, among other things, for the study of granular and composite materials (see \cite{MMZ2017} in this context), as well as \cite{S2015}, where a more general functional setting is investigated.

In the present note, we shall present the basic definitions of structured deformations and the relaxation theorems contained in \cite{CF1997,DPO1993}, and the applications contained in \cite{BMMO-expl,CMMO,MMZ2017}, which contain contributions of the author.

In Section \ref{sect:2}, we present the general mechanical and functional setting, as well as define the energies that are interesting to our problems.
In Section \ref{sect:3}, we  present the results for three specific problems.
Finally, in Section \ref{sect:4}, we outline some future directions.

%Instead of simply listing headings of different levels we recommend to
%let every heading be followed by at least a short passage of text.
%Further on please use the \LaTeX\ automatism for all your
%cross-references and citations. And please note that the first line of
%text that follows a heading is not indented, whereas the first lines of
%all subsequent paragraphs are.

\section{General setting and energies}
\label{sect:2}
In this section we recall the definition of structured deformation and also present the main results about the relaxation of non-convex energies.

Always in this note, and unless differently specified, $\Omega$ will denote a bounded open subset of the $N$-dimensional Euclidean space $\R{N}$.

\subsection{Functional setting}
\label{functionspaces}
We start by defining the function spaces which are relevant in our context.
We assume that the reader is familiar with the Lebesgue spaces $L^p(\Omega;\R d)$ ($p\in[1,\infty]$), with the space of Radon measures $\cM(\Omega;\R d)$, and with the different notions of convergence (strong, weak, and weak-*).
\begin{definition}[see \cite{AFP2000}]\label{def:BVSBV}
The space of $\R d$-valued functions of \emph{bounded variation} is defined as
%\begin{equation}\label{eq:BV}
$BV(\Omega;\R d):=\{u\in L^1(\Omega;\R d): Du\in\cM(\Omega;\R{d\times N})\}.$
%\end{equation}
The distributional derivative $Du$ of a $BV$ function is characterized by
%\begin{equation}\label{eq:distrder}
$Du=\nabla u\cL^N+[u]\otimes\nu\cH^{N-1}\res S(u)+D^cu$,
%\end{equation}
where $\nabla u$ is the approximate gradient (see \cite[Definiton 3.70 and Theorem 3.83]{AFP2000}), $[u]:=u^+-u^-$ is the jump of $u$ across the jump set $S(u)$ (which can be proved to be $(N-1)$-rectifiable), $\nu$ is the normal to $S(u)$, and $D^cu$ is the \emph{Cantor part} of the measure $Du$.
The space of \emph{special functions of bounded variation} is defined as
\begin{equation}\label{eq:SBV}
SBV(\Omega;\R d):=\{u\in BV(\Omega;\R d): D^cu=0\},
\end{equation}
the set of $BV$ functions for which the singular part of $Du$ is reduced to the jump part.
\end{definition}

We now present the two definitions of structured deformation of \cite{DPO1993} and \cite{CF1997}.
\begin{definition}[see \cite{DPO1993}]\label{def:SD-DPO}
A \emph{structured deformation} is a triple $(\kappa ,g,G)$, where $\kappa $ is a surface-like subset of $\Omega$, and the injective and piecewise differentiable map $g\colon\Omega\to\R N)$ the piecewise continuous tensor field $G:\Omega\to\R{N\times N}$ are such that
\begin{equation}
0<C<\det G(x)\leq \det \nabla g(x)\qquad\text{at each point $x\in\Omega$.}
\label{accommodation}
\end{equation}
\end{definition}
In view of Definition \ref{def:SD-DPO}, $\kappa$ describes preexisting, unopened macroscopic cracks and the map $g$ and its classical gradient $\nabla g$ describe macroscopic changes in the geometry of the body.
A geometrical interpretation of the field $G$ is provided by the following approximation theorem.
\begin{theorem}[{see \cite[Theorem 5.8]{DPO1993}}]\label{thm:approx-DPO}
For each structured deformation $(\kappa ,g,G)$ there exists a sequence of injective, piecewise smooth
deformations $f_{n}$ and a sequence of surface-like subsets $\kappa _{n}$ of
the body such that
\begin{equation}\label{eq:convergences-DPO}
g=\lim_{n\to \infty }f_{n},  \qquad G=\lim_{n\to \infty }\nabla f_{n},\qquad\text{and}\qquad \kappa =\bigcup_{n=1}^{\infty}\bigcap_{j=n}^{\infty }\kappa_{j}.
\end{equation}
\end{theorem}
The limits in \eqref{eq:convergences-DPO} are taken in the sense of $L^{\infty}$convergence. 
From \eqref{eq:convergences-DPO} we see that $G$ captures the effects at the macroscopic level of smooth geometrical changes at submacroscopic levels; $G$ is usually referred to as the \emph{deformation without disarrangements}.

The space $SBV(\Omega;\R d)$ defined in \eqref{eq:SBV} formalizes the notion of discontinuity of a function $u$ with the introduction of the jump set $S(u)$.
In this way, one can think of linking the role of $\kappa$ with that of $g$ in Definition \ref{def:SD-DPO}, by considering pairs made of a deformation $g$, also carrying the information about the cracks, and a matrix-valued field $G$.
\begin{definition}[see \cite{CF1997}]\label{def:SD-CF}
The space of structured deformations is defined as
\begin{equation}\label{eq:SD-CF}
SD(\Omega):=\{(g,G): g\in SBV(\Omega;\R d), G\in L^1(\Omega;\R{d\times N})\}.
\end{equation}
\end{definition}
In view of Definition \ref{def:SD-CF}, Theorem \ref{thm:approx-DPO} has the following counterpart.
\begin{theorem}[{see \cite[Theorem 2.12]{CF1997}}]\label{thm:approx-CF}
Let $(g,G)\in SD(\Omega)$.
Then there exist $u_n\in SBV(\Omega;\R d)$ such that
\begin{equation}\label{eq:convergences-CF}
u_n\to g\;\text{in $L^1(\Omega;\R d)$,}\qquad \nabla u_n\wsto G\;\text{in $\cM(\Omega;\R{d\times N})$.}
\end{equation}
\end{theorem}
The proof of Theorem \ref{thm:approx-CF} can be achieved by combining Alberti's theorem \cite[Theorem 3]{A1991} and an approximation result in $BV$ by piecewise constant functions (see, \emph{e.g.}, \cite[Lemma 2.9]{CF1997}).
It is worth mentioning that an alternative proof of a weakerversion of the approximation theorem, without using Alberti's theorem, is proposed in \cite[Theorem 7.1]{S2015}. 
We want to stress that the convergences in \eqref{eq:convergences-CF} allow for the limit of a sequence of functions and the limit of its gradients to be unrelated, and in fact there is no \emph{a priori} relationship between $g$ and $G$.
Even more importantly, the \emph{disarrangements tensor} defined by the difference
\begin{equation}\label{eq:M}
M:=\nabla g-G
\end{equation}
has a fundamental geometrical meaning: it captures, in the limit as $n\to\infty$, the volume density of separations and slips between pieces of the body, and how this is determined by the interfacial discontinuities of the approximating deformations $u_n$ in Theorem \ref{thm:approx-CF} (or $f_n$ in Theorem \ref{thm:approx-DPO}).
The tensor $M$ can also be interpreted as a measure of how non classical a deformation is: $M=0$ corresponds to $\nabla g=G$, and therefore the convergences in \eqref{eq:convergences-DPO} and \eqref{eq:convergences-CF} hold in a stronger sense.
In this regard, from \eqref{eq:M} and \eqref{eq:convergences-DPO} one trivially obtains $M=\nabla \big(\lim_{n\to\infty f_n}\big)-\lim_{n\to\infty} \nabla f_n$, which can be considered a quantitative measure of the lack of commutativity of the classical gradient and the limit operator in the $L^\infty$ convergence.

The tensor $M$ and its derivatives, such as $\curl M$, are fundamental quantities to describe the presence of defects such as dislocations, and how they intervene in the response of solids (see the discussion in \cite[Section 1.1]{BMMO-expl}).

We close this discussion by presenting two examples (see \cite{CF1997,DPO1993})
\begin{example}[The broken ramp]\label{ex:BR}
Let $N=1$, $\Omega=(0,1)$, $\kappa=\emptyset$, $g(x)=2x$, and $G(x)=1$.
An approximating sequence is given by $f_n(x):=x+\frac kn$, for $\frac kn\leq x<\frac{k+1}n$ and $k=0,\ldots,n-1$.
Indeed, $f_n(x)\to2x$ as $n\to\infty$, and $\nabla f_n(x)=1$ for every $n$, so that \eqref{eq:convergences-DPO} and \eqref{eq:convergences-CF} are satisfied; yet, the distributional derivative is given by $Df_n=1+\sum_{k=1}^{n-1}\frac1k\delta_{k/n}$ and it shows the emergence of jumps discontinuities, which become smaller and smaller in magnitude and diffuse in the whole domain as $n\to\infty$, see Figure \ref{fig:one}.
\begin{figure}[h]
\begin{center}
\includegraphics[scale=.7]{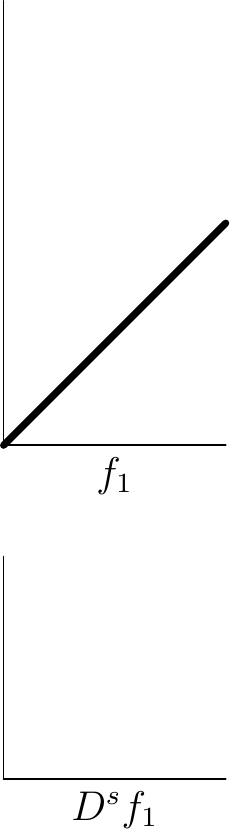}$\qquad$ 
\includegraphics[scale=.7]{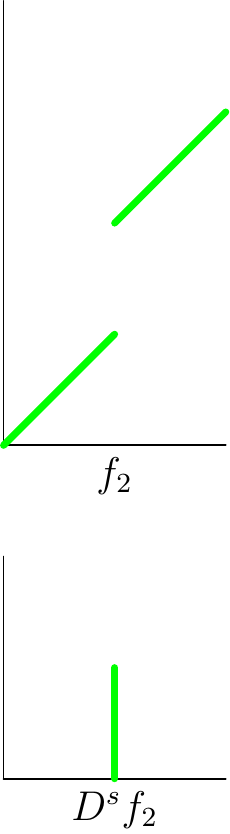}$\qquad$ 
\includegraphics[scale=.7]{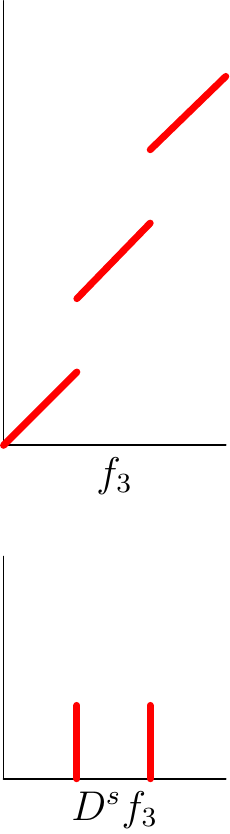}$\qquad$ 
\includegraphics[scale=.7]{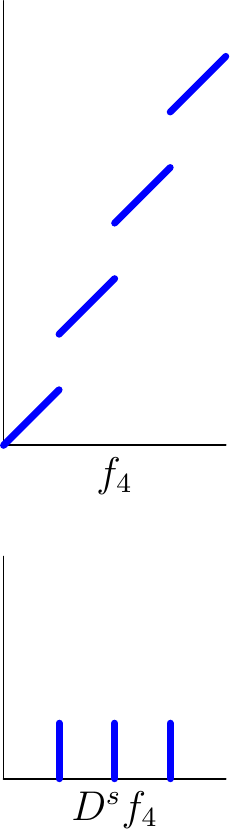}%$\qquad\qquad$ 
\caption{The first four steps of the approximation in Example \ref{ex:BR}. The top row depicts the approximants $f_1,\ldots,f_4$; the bottom row depicts the singular parts $D^sf_1,\ldots,D^sf_4$.}
\label{fig:one}
\end{center}
\end{figure}
\end{example}
\begin{example}[The deck of cards]\label{ex:DC}
In this example in $N=3$ dimensions, a simple shear of a cubic domain $\Omega=(0,1)^3$ is considered.
Let $\kappa=\emptyset$, $g(x)=(x_1+x_3,x_2,x_3)$, and $G(x)=\mathbb I$.
An approximating sequence is given by $f_n(x):=(x_1+\frac kn,x_2,x_3\big)$, for $\frac kn\leq x_3<\frac{k+1}n$ and $k=0,\ldots,n-1$.
Then, $f_n(x)\to g(x)$ as $n\to\infty$, and $\nabla f_n(x)=\mathbb I$ for every $n$, so that \eqref{eq:convergences-DPO} and \eqref{eq:convergences-CF} are once again satisfied; here, $Df_n=\mathbb I+\sum_{k=1}^{n-1}\frac1k\delta_{k/n}(x_3)e_1\otimes e_3$, and also in this case the jumps diffuse to the bulk.
\end{example}

The effect that the singular part of $Du_n$ diffuses to the bulk is crucial to understanding structured deformations and the cell formulas in the relaxation theorems that we will present later.
What is happening is that, as $n\to\infty$, the singular sets diffuse ($\cH^{N-1}(S(u_n))\to+\infty$) in a controlled way, which is measured by the boundedness of the total variation of the singular part of the measure, $|D^su_n|(\Omega)<+\infty$.

\subsection{Energies and relaxation}
\label{energies}
We present here some general energy functionals for deformations of solids, as well as the relaxation theorems to define functionals on the space $SD(\Omega)$ introduced in \eqref{eq:SD-CF}.
Energy minimization is a convenient way to find equilibrium configurations of a body undergoing internal stresses and external loadings.
In the general case, the energy densities are not convex, or quasi-convex, so that the energy landscape can show many local minima. 
In calculus of variations, this situation is resolved by computing the (quasi-)convex envelope of the energy densities, therefore obtaining the relaxed functional, that is the largest lower semicontinuous functional below the given one.
We refer the reader to treatises on calculus of variations for the relationship between lower semicontinuity and convexity (see, \emph{e.g.}, \cite{FL2007}).

Let now $u\in SBV(\Omega;\R d)$ and consider an energy functional defined by
\begin{equation}\label{eq:en_fun}
E(u) := \int_{\Omega} W(\nabla u(x))\,\de x+\int_{S(u)\cap\Omega}\Psi([u],\nu(u))\,\de\cH^{N-1},
\end{equation}
where the desities $W\colon\R{d\times N}\to[0,+\infty)$ and $\Psi\colon\R d\times\S{N-1}\to[0,+\infty)$ satisfy: %the following hypotheses
\begin{enumerate}
\item[$(H_1)$] $\exists\,C>0$ such that $\forall\,A,B\in\R{d\times N}$:
$|W(A)-W(B)|\leq C|A-B|(1+|A|^{p-1}+|B|^{p-1})$
for some $p\in(1,+\infty)$;
\item[$(H_2)$] $\exists\,c_1,C_1>0$ such that $\forall\,(\lambda,\nu)\in\R d\times\S{N-1}$:
$c_1|\lambda| \leq \Psi(\lambda, \nu) \leq C_1|\lambda|$;
\item[$(H_3)$] (positive homogenity of degree $1$) $\forall\,\lambda\in\R d$, $t>0$:
$\Psi(t\lambda,\nu)=t\Psi(\lambda,\nu)$;
\item[$(H_4)$]  (subadditivity) $\forall\,\lambda_1,\lambda_2\in\R d$:
$\Psi(\lambda_1+\lambda_2,\nu)\leq\Psi(\lambda_1,\nu)+\Psi(\lambda_2,\nu)$.
\end{enumerate}
Given a structured deformation $(g,G)\in SD(\Omega)$, the relaxation of the energy \eqref{eq:en_fun} is %defined as
\begin{equation}\label{eq:rel_en}
\begin{split}
I(g,G):=\inf\bigg\{ & \liminf_{n\to\infty} E(u_n):\, u_n\in SBV(\Omega;\R d),\; u_n\to g \text{ in } L^1(\Omega;\R d), \\
& \nabla u_n\wsto G \text{ in } \cM(\Omega;\R{d\times N}),\; \sup_n \lVert\nabla u_n\rVert_{L^p}<\infty\bigg\}.
\end{split}
\end{equation}
\begin{remark}\label{rem:densityhp}
We notice that the lack of coercivity of the bulk density $W$ is compensated by the boundedness constraint on the approximating functions $u_n$ requested in \eqref{eq:rel_en}.
This gives the freedom to treat more general densities $W$ satisfying $(H_1)$.
A strategy to circumvent the lack of coercivity is proposed in \cite[proof of Proposition 2.22, Step 2]{CF1997}: one first adds a fictitious coercivity $W^\eps(\cdot):=W(\cdot)+\eps|\cdot|^p$, and then lets $\eps\to0$.
The requirement in $(H_2)$ of coercivity for $\Psi$ in the first variable can be removed \cite[Remark 3.3]{CF1997} by modifying the requirements in $(H_2)$ and $(H_3)$ and the boundedness requirement in \eqref{eq:rel_en}.
\end{remark}
Let us denote by $Q:=(-\frac12,\frac12)^N$ the unit cube in $\R N$ and by $Q_\nu$ the rotated one such that two faces are perpendicular to the vector $\nu\in\S{N-1}$.
The main result regarding the relaxed energy in \eqref{eq:rel_en} is the following.
\begin{theorem}[{see \cite[Theorem 2.17]{CF1997}}]\label{thm:rel-CF}
Let $(g,G)\in\ SD(\Omega)$ with $G\in L^p(\Omega;\R{d\times N})$, and let the bulk and surface energy densities $W$ and $\Psi$ satisfy $(H_1)$-$(H_4)$.
Then
\begin{equation}\label{eq:int_rep-CF}
I(g,G)=\int_\Omega H(\nabla g,G)\,\de x+\int_{S(g)\cap\Omega} h([g],\nu(g))\,\de\cH^{N-1},
\end{equation}
where, for $A,B\in\R{d\times N}$, $\lambda\in\R d$, $\nu\in\S{N-1}$
\begin{equation}\label{eq:H-CF}
\begin{split}
H(A,B):= \inf\bigg\{ & \int_Q W(\nabla u)\,\de x+\int_{S(u)\cap Q} \Psi([u],\nu(u))\,\de\cH^{N-1}: \\
& u\in SBV(Q,\R d),\; u|_{\partial Q}=Ax,\; |\nabla u|\in L^p(\Omega),\; \int_Q \nabla u\,\de x=B\bigg\},
\end{split}
\end{equation}
%and for $\lambda\in\R d$, $\nu\in\S{N-1}$
\begin{equation}\label{eq:h-CF}
\begin{split}
h(\lambda,\nu):= \inf\bigg\{ & \int_{S(u)\cap Q_\nu} \Psi([u],\nu(u))\,\de\cH^{N-1}: \\
& u\in SBV(Q_\nu;\R d),\; u|_{\partial Q_\nu}=u_{\lambda,\nu}(x),\;\nabla u=0\;\cL^N\text{-a.e.}\bigg\},
\end{split}
\end{equation}
with
%\begin{equation}\label{ulambdanu}
$u_{\lambda,\nu}(x):=
\begin{cases}
0 & \text{if $-\frac12\leq x\cdot\nu<0$,} \\
\lambda & \text{if $0\leq x\cdot\nu<\frac12$.}
\end{cases}$
%\end{equation}
\end{theorem}
The relaxation described in Theorem \ref{thm:rel-CF} is based on the so-called \emph{blow-up method} \cite{FM1993}: the energy densities $H$ and $h$ are obtained as the energetically optimal ones realizing an affine transformation (for the bulk part) with prescribed boundary values and gradient, and a jump (for the surface part) with prescribed jump height.
This is essentially the meaning of the conditions in the infimization problems \eqref{eq:H-CF} and \eqref{eq:h-CF}.

Results analogous to that stated in Theorem \ref{thm:rel-CF} are available in the literature: the reader is referred to \cite{CF1997} for the case $p=1$, and to \cite{BMMO,OP2000,P2004} for second-order structured deformations and for different types of initial energies \eqref{eq:en_fun}.

\section{Three problems}
\label{sect:3}
%In this section, we collect three problems that are useful in applications. 
\subsection{Explicit formulas for purely interfacial energies}
\label{sect:3.1}
In the case $d=N$, we investigate an energy functional of the type \eqref{eq:en_fun} with $W=0$, that is, where the energy is purely interfacial.
By applying the blow-up method, we obtain a representation formula of the relaxed energy \eqref{eq:rel_en} for particular choices of initial surface energy densities, namely for $E^{|\cdot|}(u):=\int_{S(u)\cap\Omega} |[u]\cdot\nu(u)|\,\de\cH^{N-1}$ and $E^{\pm}(u):=\int_{S(u)\cap\Omega} ([u]\cdot\nu(u))^\pm\,\de\cH^{N-1}$.
Given $(g,G)\in SD(\Omega)$, denote by $\cV^{|\cdot|}(g,G)$ and $\cV^\pm(g,G)$ the relaxed energies obtained via Theorem \ref{thm:rel-CF}.
%$\Psi^{\vert \cdot \vert}(\lambda ,\nu ):=|\lambda\cdot\nu|$, and $\Psi^{\pm }(\lambda,\nu ):= (\lambda\cdot\nu)^\pm$ (positive and negative part).
%The results contained in \cite{BMMO-expl} recover those from \cite{OP2015} for the same explicit energy densities, and provide the equivalence of two minimum problems, the first of which is the adaptation of \eqref{eq:H-CF} to the present case, and the second of which is a differen relaxation process proposed in \cite{BMS2012}.

\begin{theorem}[see \cite{BMMO-expl,OP2015,S}]\label{thm:expl}
The initial disarrangement densities
%\begin{equation}\label{eq:Psi-all}
$\Psi ^{\vert \cdot \vert}(\lambda ,\nu ):=\vert \lambda \cdot \nu\vert$ and $\Psi ^{\pm }(\lambda ,\nu ):=(\lambda \cdot \nu )^{\pm}$
%\end{equation}%
%satisfy the hypotheses $(H_2)$-$(H_4)$ and 
have relaxed
disarrangement densities given by
%\begin{equation}
$H^{\lvert \cdot\rvert}(A,B)=\vert \tr(A-B)\vert$ , %\quad 
$h^{\lvert \cdot \rvert }(\lambda ,\nu )=\Psi ^{\lvert \cdot \rvert }(\lambda ,\nu )$,
%\end{equation}%
and
%\begin{equation}
$H^{\pm }(A,B)=(tr(A-B))^{\pm}$, $h^{\pm }(\lambda ,\nu )=\Psi ^{\pm}(\lambda ,\nu )$.
%\label{plus or minus relaxed densities}
%\end{equation}
\end{theorem}
The proof presented in \cite{BMMO-expl} that $H^{\lvert \cdot\rvert}(A,B)=\vert \tr(A-B)\vert$ relies on the following chain of inequalities
\begin{equation}\label{eq:chain}
\begin{split}
\lvert \tr(A-B)\rvert  \leq \inf \bigg\{ & \int_{S(u)\cap Q}\lvert [u](x)\cdot \nu(u)(x)\rvert\,\de\cH^{N-1}(x):\;u\in SBV(Q;\R N),  \\
& u|_{\partial Q}=Ax,\;\nabla u\in L^{p}(Q),\;\int_{Q}\nabla u\,\de x=B\bigg\}   \\
\leq \inf \bigg\{& \int_{S(u)\cap Q}\lvert [u](x)\cdot \nu(u)(x)\rvert\,\de\cH^{N-1}(x):\;u\in SBV(Q;\R N), \\
& u|_{\partial Q}=0,\;\nabla u=B-A\text{ a.e. in } Q\bigg\}\leq|\tr(A-B)|.
\end{split}
\end{equation}
The first inequality comes from \eqref{eq:H-CF}, the second one holds because we are restricting the set of admissible functions in the minimization, the third one is proved in \cite{BMMO-expl}.
It is worth noticing that the second minimization problem in \eqref{eq:chain} is the one proposed in \cite{BMS2012} in a different context for relaxation; the proof of the third inequality in \eqref{eq:chain} is achieved by an explicit construction and on the notion of \emph{isotropic vectors} %discussed in 
\cite{CL1998}.

Other than recovering the explicit results of \cite{OP2015} with a shorter proof, and other additional explicit formulas, \eqref{eq:chain} shows the equivalence of two minimum problems.
The recent paper \cite{S} contains results on the general form of the relaxation of purely interfacial energies, which extends the previous results.

Another consequence is the relationship between the relaxed energied for the three densities considered in Theorem \ref{thm:expl}, namely (recalling \eqref{eq:M})
$$\cV^\pm(g,G)=\frac12\cV^{|\cdot|}(g,G)\pm\frac12\int_{\Omega}\tr M(x)\,\de x.$$

\subsection{Optimal design}
\label{sect:3.2}
In the context of optimal design, we consider a two-component fractured medium with prescribed macroscopic strain.
The initial energy functional of type \eqref{eq:en_fun} is tailored to account for the fine structure of the material and the relaxed energy densities are obtained by the interplay between the optimization of sharp interfaces and the diffusion of microscopic cracks.

Let $\chi\in BV(\Omega;\{0,1\})$ be the characteristic function of a set of finite perimeter (see \cite{AFP2000}) describing one constituent of the material.
Given a deformation $u\in SBV(\Omega;\R d)$, consider the initial energy
\begin{equation}\label{energyodsd}
\begin{split}
E(\chi, & u):= \int_\Omega ((1-\chi)W^0(\nabla u)+\chi W^1(\nabla u) )\,\de x \\
+ &  \int_{\{ \chi = 0\}\cap S(u)\cap\Omega}\Psi^0_1([u],\nu(u))\,\de\cH^{N-1}+ \int_{\{ \chi = 1\}\cap S(u)\cap\Omega} \Psi^1_1([u],\nu(u))\,\de\cH^{N-1} \\
+ &  \int_{S(\chi)\cap S(u)\cap\Omega} \Psi_2(\chi^+,\chi^-,u^+,u^-,\nu(u))\,\de\cH^{N-1} + |D\chi|(\Omega),
\end{split}
\end{equation}
which features the contributions of bulk and surface energy densities $W^i$ and $\Psi_1^i$, $i=0,1$, for the deformation $u$, as well as a surface energy density $\Psi_2$ on the superposition of the singular sets of $u$ and $\chi$, and a perimeter penalization term pushing for smaller interfaces between the two constituents.

We assume that the energy densities $W^i$ satisfy $(H_1)$ and are coercive, that $\Psi_1^i$ satisfy $(H_2)$-$(H_4)$, and that $\Psi_2$ satisfies
\begin{enumerate}
\item[$(H_5)$] $\exists\,C>0$ such that $\forall\,a,b\in\{0,1\}$, $c, d\in\R{d}$, $\nu \in \S{N-1}$:\\
$0 \leq \Psi_2(a,b,c,d, \nu) \leq C(1+ |a - b|+|c-d|)$;
\item[$(H_6)$] $\forall\,a, b\in \{0,1\}$, $c,d \in \R{d}$, $\nu \in \S{N-1}$:
$\Psi_2(a, b,c,d, \nu) = \Psi_2(b,a,  d, c,-\nu)$;
\item[$(H_7)$] $\exists\,C>0$ such that $\forall\,a,b, \in \{0,1\}$, $c_i, d_i \in \R{d}$, $i=1, 2$, $\nu \in \S{N-1}$,\\
$|\Psi_2(a,b, c_1, d_1, \nu)- \Psi_2(a,b, c_2, d_2, \nu)| \leq %C \big||c_1-c_2|- |d_1-d_2|\big|\leq 
C|(c_1-d_1)-(c_2-d_2)|$;
\item[$(H_8)$] $\forall\,a\in \{0,1\}$, $c\in \R{d}$: $\Psi_2(\cdot, \cdot, c, c, \cdot) = \Psi_2(a,a,\cdot,\cdot, \cdot)=0$.
\end{enumerate}
Given $(\chi,g,G)\in BV(\Omega;\{0,1\})\times SD(\Omega)$, the relaxed energy is defined by
\begin{equation}\label{eq:rel_fractured}
\begin{split}
I(\chi, & u,G):=  \inf\bigg\{\liminf_{n\to\infty} E(\chi_n,u_n): \chi_n\in BV(\Omega; \{0,1\}),\; u_n\in SBV(\Omega;\R d), \\
& \chi_n\wsto \chi \text{ in } BV(\Omega; \{0,1\}),\; u_n\to u \text{ in } L^1(\Omega;\R d) ,\nabla u_n\wto G \text{ in }L^p(\Omega;\R {d\times N}) \; \bigg\},
\end{split}
\end{equation}
and the relaxation theorem states the following.
\begin{theorem}[{\cite[Theorem 3.3]{MMZ2017}}]\label{thm:rel_fractured}
Let $(\chi,g,G)\in BV(\Omega; \{0,1\})\times SD(\Omega)$ and let the bulk and surface energy densities $W^i$, $\Psi_1^i$ ($i=0,1$), and $\Psi_2$ satisfy $(H_1)$-$(H_8)$.
Then
\begin{equation}\label{intrep}
I(\chi,g,G)= \int_\Omega H ( \chi, \nabla g, G)\, \de x + \int_{S(\chi,g)\cap\Omega} h( \chi^+, \chi^-, g^+, g^-, \nu)\, \de \cH^{N-1}
\end{equation}
where, for $i,a,b\in \{0,1\}$, $A, B \in \R{d{\times}N}$, $c, d \in \R{d}$, and $\nu \in \S{N-1}$,
\begin{equation}\label{Hsd}
\begin{split}
H (i, A, B):=\inf\bigg\{ & \int_Q W^i(\nabla u)\, \de x + \int_{S(u)\cap Q} \Psi^i_1( [u], \nu(u)) \, \de \cH^{N-1}: \\
&  u \in SBV( Q; \R{d}),\; |\nabla u| \in L^p(Q), \; u|_{\partial Q} = Ax, \; \int_Q \nabla u\,\de x = B \bigg\},
\end{split}
\end{equation}
%and, for $a, b \in \{0,1\}$, $c, d \in \R{d}$, $\nu \in \cS^{N-1}$,
\begin{equation}\label{gamma} 
\begin{split}
h( a, b, c, d, \nu) &:= \inf \bigg\{ \int_{S(\chi)\cap S(u)\cap Q_\nu} \Psi_2(\chi^+, \chi^-, u^+, u^-, \nu(u))\, \de \cH^{N-1} +|D \chi|(Q_\nu) \\
& + \int_{\{ \chi = 0\} \cap S(u)\cap Q_\nu} \Psi^0_1([u], \nu(u))\, \de \cH^{N-1} \\
& + \int_{\{ \chi = 1\} \cap S(u)\cap Q_\nu} \Psi_1^1([u], \nu(u)) \, \de \cH^{N-1}:  (\chi,u) \in \cA(a,b, c, d,\nu)\bigg\},
\end{split}
\end{equation}
where $\cA(a,b,c, d, \nu) := \{(\chi, u) \in BV(Q_\nu; \{0,1\}){\times}SBV(Q_\nu; \R{d}): \chi|_{\partial Q_\nu} =\chi_{a,b,\nu},$ $u|_{\partial Q_\nu} = u_{c,d, \nu}, \;\nabla u = 0 \;\mathcal{L}^N\text{-a.e.}\}$
with 
$$ \chi_{a,b, \nu}(x) :=\begin{cases} 
	a & \text{if } x\cdot \nu > 0, \\
	b & \text{if } x\cdot \nu \leq 0,
	\end{cases}
\qquad\text{and}\quad 
u_{c,d, \nu}(x) := \begin{cases} 
	c & \text{if } x\cdot \nu > 0, \\
	d & \text{if } x\cdot \nu \leq 0.
	\end{cases}
$$
\end{theorem}

\subsection{Dimension reduction}
\label{sect:3.3}
Dimension reduction is a technique to study thin objects in $2$d starting from a fully $3$d model and energetics.
A small thickness parameter $\eps>0$ is involved in the $3$d energy and sent to zero to obtain the $2$d energy.
We consider a bounded open set $\omega\subset\mathbb{R}^{2}$ and define $\Omega_{\eps}:=\omega\times(-\frac{\eps}{2},\frac{\eps}{2})$; for $u\in SBV(\Omega_\eps;\R3)$ let
\begin{equation}\label{eq:E3d}
E_{\eps}(u):=\int_{\Omega_{\eps}}W_{3d}(\nabla u(x))\,\de x+\int_{S(u)\cap\Omega_{\eps}} \Psi_{3d}([u](x),\nu(u)(x))\,\de\cH^{2},
\end{equation}
where the densities $W_{3d}$ and $\Psi_{3d}$ satisfy assumptions $(H_1)$-$(H_4)$.

\begin{figure}[h]
\begin{center}
\begin{tikzpicture}
  \matrix (m) [matrix of math nodes,row sep=2em,column sep=3em,minimum width=1em]
  {
    & W_{3d}, h_{3d} & \\
    W_{3d,2d}, h_{3d,2d} & & W_{3d,SD}, h_{3d,SD} \\
    W_{3d,2d,SD}, h_{3d,2d,SD} & & W_{3d,SD,2d}, h_{3d,SD,2d} \\
    & W_1,\Gamma_1 & \\
    };
  \path[-stealth]
    (m-1-2) edge [double] node [left] {$\mathrm{DR}\quad$} (m-2-1)
            edge [double] node [right] {$\quad\mathrm{SD}$} (m-2-3)
    (m-2-1) edge [double] node [left] {$\mathrm{SD}$} (m-3-1)
    (m-2-3) edge [double] node [right] {$\mathrm{DR}$} (m-3-3)
    (m-1-2) edge node [left] {\cite{MS2014}} (m-4-2);
\end{tikzpicture}
\caption{Energy densities for the paths for dimension reduction and structured deformations in \cite{CMMO}.}
\label{figtwo}
\end{center}
\end{figure}
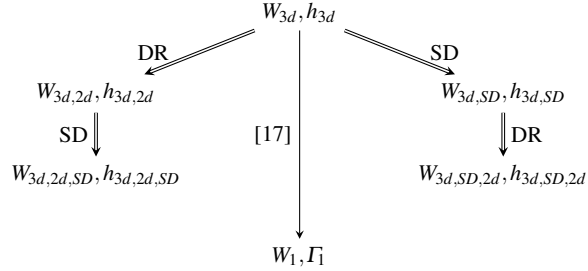

In order to couple dimension reduction and structured deformations, we shall perform the two relaxation processes one after the other in two different orders, as depicted in the left- and right-hand side paths in Fig. \ref{figtwo}.
The dimension reduction legs DR are obtained by the change of variables $x_3\to x_3/\eps$, so that the functional will be defined on the volume $\Omega:=\Omega_1=\omega\times(-\frac12,\frac12)$.
After rescaling the functional by dividing it by $\eps$, the limit as $\eps\to0$ is computed.
The structured deformation legs SD are obtained by applying the relaxation Theorem \ref{thm:rel-CF}.
It is important to point out that, in the dimension reduction process, there is the emergence of a vector $\cl d:\omega\to\R3$ that keeps memory of the behavior of the deformation in the $x_3$-direction.
The role played by $\cl d$ is encoded via the constraint $\int_{-1/2}^{1/2} \frac{\nabla_3 u_n}{\eps_n}\,\de x_3\wto\cl d$ in the relaxation process (see \cite{CMMO}).

A comparison of the relaxed energy densities is interesting, to see whether or not the two different ways provide the same result.
In this respect, partial results are available in the special case of the purely interfacial initial energy densities in Theorem \ref{thm:expl}.
%purely interfacial energies $W_{3d}=0$ (where also the results presented in Section \ref{sect:3.1} can be used).
Moreover, in this case, another relaxation procedure is available in the literature: in \cite{MS2014} a relaxation that simultaneously defines a $2$d energy on structured deformations is studied.
We prove \cite[Section 5]{CMMO} that the left- and right-hand paths provide the same relaxed energy densities, whereas those computed using the central path in Fig. \ref{figtwo} are lower.

\section{Conclusions and outlook}
\label{sect:4}
In this brief note, we have presented a general overview of structured deformations and their versatility by describing three applications which are very different from one another. 
In Section \ref{sect:3.1} we obtained explicit formulas and the equivalence of two minimum problems \cite{BMMO-expl}; in Section \ref{sect:3.2} we coupled structured deformations to study an optimal design problem, and performed a relaxation in the joint variables $\chi$ and $u$ \cite{MMZ2017}; finally, in Section \ref{sect:3.3} we coupled structured deformations with dimension reduction, showing two possible routes to obtain the final relaxed energies, and comparing the results for a specific choice of initial energy densities \cite{CMMO}.

The results in Section \ref{sect:3.1} can be extended to more general situations and call for the search of explicit formulas for other initial energy densities: these will be certainly useful in the applications.
Those in Section \ref{sect:3.2} are a first step towards the far-reaching goal of incorporating elements of plasticity in optimal design of composite media.
Section \ref{sect:3.3} presents challenges at a more theoretical level, regarding the commutativity of different relaxation processes and the comparison of minimum problems for the same functional, but different sets of admissible functions.

The versatility of structured deformations makes them suitable to model physical systems such as composite or granular media, defective materials, and biological membranes with the aim of shedding new light on the mechanics of these systems.

\begin{acknowledgement}
The author acknowledges partial support for this research from the following grants: FCT$\_$UTA/CMU/MAT/0005/2009 of the Funda\c{c}\~{a}o para a Ci\^encia e a Tecnologia through the Carnegie Mellon Portugal Program; ERC Advanced grant \emph{Quasistatic and Dynamic Evolution Problems in Plasticity and Fracture} (Grant agreement 290888); INdAM-GNAMPA project 2015 \emph{Fenomeni Critici nella Meccanica dei Materiali: un Approccio Variazionale}; ERC Starting grant \emph{High-Dimensional Sparse Optimal Control} (Grant agreement 306274).
The author is a member of the %Gruppo Nazionale per l'Analisi Matematica, la Probabilit\`a e le loro Applicazioni (
GNAMPA group %) 
of %the Istituto Nazionale di Alta Matematica (
INdAM. %).
The author is grateful to J. Matias and D. R. Owen for valuable suggestions in writing this note.
%is thankful to the Department of Mathematical Sciences at Carnegie Mellon University (Pittsburgh, PA, USA), the Departamento de Matem\'atica at Instituto Superior T\'ecnico (Lisbon, Portugal), the Mathematics Area at SISSA (Trieste, Italy), and the Fakult\"at f\"ur Mathematik at Technische Universit\"at M\"unchen (Munich, Germany), where the research leading to these results has been conducted with the partial support of the grants 
\end{acknowledgement}
%

%\input{referenc}
%% REFERENCES %%

% Silhavy -- quello under review. (altri?)
% De Giorgi - Letta -- criterion
% Deseri - Owen -- qualcosa?

\begin{thebibliography}{99}
%\bibitem{AF2015} A. Acharya and C. Fressengeas. %: \emph{Continuum mechanics of the interaction of phase boundaries and dislocations in solids}. 
%In \emph{Differential Geometry and Continuum Mechanics}, G. Q Chen, M. Grinfeld, R.J. Knops eds., Springer Proceedings in Mathematics \& Statistics \textbf{137} (2015), 125-168.

%\bibitem{AD2015} V. Agrawal and K. Dayal. %: \emph{Dynamic Phase-field Model for Structural Transformations and Twinning: Regularized Interfaces with Transparent Prescription of Complex Kinetics and Nucleation. Part I: Formulation and One-Dimensional Characterization}. 
%J. Mech. Phys. Solids \textbf{85} (2015), 270-290.

\bibitem{A1991} G. Alberti. %: \emph{A Lusin type theorem for gradients}, 
J. Funct. Anal., \textbf{100} (1991), 110-118.

\bibitem{AFP2000} L. Ambrosio, N. Fusco, and D. Pallara. %: \emph{Functions of Bounded Variation and Free Discontinuity Problems}, 
OUP, 2000. ISBN: 9780198502456.

%\bibitem{BMS2011} M. Ba\'ia, J. Matias and P. M. Santos: \emph{A survey on structured deformations}. Sao Paulo J. Math. Sci., \textbf{5} (2011), 185-201.

\bibitem{BMS2012} M. Ba\'ia, J. Matias, and P. M. Santos. %: \emph{A relaxation result in the framework of structured deformations in a bounded variation setting}. 
Proc. Royal Soc. Edinburgh, \textbf{142A} (2012), 239-271.

\bibitem{BMMO-expl} A. C. Barroso, J. Matias, M. Morandotti, and D. R. Owen. %: \emph{Explicit Formulas for Relaxed Disarrangement Densities Arising from Structured Deformations}. 
Preprint SISSA 37/2015/MATE. %arXiv:1508.06908. 
Math. Mech. Complex Syst., \emph{to appear}.

\bibitem{BMMO} A. C. Barroso, J. Matias, M. Morandotti, and D. R. Owen. %: \emph{Second-order structured deformations: relaxation, integral representation and applications}. 
Preprint SISSA 37/MATE/2016. %arXiv:1607.02311. \emph{Submitted}.

\bibitem{CMMO} G. Carita, J. Matias, M. Morandotti, and D. R. Owen: \emph{Dimension reduction in the context of first-order structured deformations}. \emph{In preparation}

\bibitem{CDPFO1999} R. Choksi, G. Del Piero, I. Fonseca, and D. R. Owen. %: \emph{Structured deformations as energy minimizers in models of fracture and hysteresis}. 
Math. Mech. Solids., \textbf{4} (1999), 321-356.

\bibitem{CF1997} R. Choksi and I. Fonseca. %: \emph{Bulk and interfacial energy densities for structured deformations of continua}. 
Arch. Rational Mech. Anal., \textbf{138} (1997), 37-103.

\bibitem{CL1998} N. Ciblak and H. Lipkin. %: \emph{Orthonormal Isotropic Vector Bases}. 
Proceedings of DETC'98, September 13-16, Atlanta, Georgia.

\bibitem{DPO1993} G. Del Piero and D. R. Owen. %: \emph{Structured deformations of continua}. 
Arch. Rational Mech. Anal., \textbf{124} (1993), 99-155.

\bibitem{DPT2001} G. Del Piero and L. Truskinovsky. %: \emph{Macro- and micro-cracking in one-dimensional elasticity}. 
Int. J. of Solids and Structures, \textbf{38} (2001), 1135-1148.

\bibitem{DO2003} L. Deseri and D. R. Owen. %: \emph{Toward a field theory for elastic bodies undergoing disarrangements}. 
J.\@ Elast., \textbf{70} (2003), 197-236.

\bibitem{DO2015} L. Deseri and D. R. Owen. %: \emph{Stable disarrangement phases arising from expansion/contraction or from simple shearing of a model granular medium}. 
Int. J. Eng. Sci. \textbf{96} (2015), 111-140.

\bibitem{FL2007} I. Fonseca and G. Leoni. %: \emph{Modern methods in the calculus of variations: $L^p$ spaces}. Springer Monographs in Mathematics. 
Springer, New York, 2007. ISBN 978-0-387-69006-3.

\bibitem{FM1993} I. Fonseca and S. M\"uller. %: \emph{Relaxation of quasi-convex functionals in $BV(\Omega;\R p)$ for integrands $f(x,u,\nabla u)$}. 
Arch. Rational Mech. Anal. {\bf 123}, (1993), 1-49.

%\bibitem{JH2000} R. D. James and K. F. Hane. %: \emph{Martensitic transformations and shape-memorymaterials}. 
%Acta Mater. \textbf{48} (2000), 197-222.

%\bibitem{M2007} J. Matias: \emph{Differential inclusions in $SBV_{0}(\Omega )$ and applications to the calculus of variations}. J. Convex Analysis \textbf{14} (2007) No. 3, 465-477.

\bibitem{MMZ2017} J. Matias, M. Morandotti, and E. Zappale. %: \emph{Optimal Design of Fractured Media with Prescribed Macroscopic Strain}. 
J. Math. Anal. Appl., \textbf{449} (2017), 1094-1132.
%\DOI{10.1016/j.jmaa.2016.12.043}.

\bibitem{MS2014} J. Matias and P. M. Santos. %: \emph{A dimension-reduction result in the context of structured deformations}, 
Appl. Math. Optim., \textbf{69} (2014), 459-485.

\bibitem{O} D. R. Owen. %\emph{Elasticity with Gradient-Disarrangements: a Multiscale Geometrical Perspective for Strain-Gradient Theories of Elasticty and of Plasticity}. %
J. Elasticity, DOI: 10.1007/s10659-016-9599-9.

\bibitem{OP2000} D. R. Owen and R. Paroni. %: \emph{Second-order structured deformations}. 
Arch. Rational Mech. Anal., \textbf{155} (2000), 215-235.

\bibitem{OP2015} D. R. Owen and R. Paroni. %: \emph{Optimal flux densities for linear mappings and the multiscale geometry of structured deformations}. 
Arch. Rational Mech. Anal., \textbf{218} (2015), 1633-1652.

\bibitem{P2004} R. Paroni. %: \emph{Second-order structured deformations: approximation theorems and energetics}, 
In %\emph{Multiscale modeling in continuum mechanics and structured deformations}, 
G. Del Piero and D. R. Owen eds., CISM \textbf{447}, Springer (2004).

\bibitem{S2015} M. \v{S}ilhav\'y. %: \emph{On the approximation theorem for structured deformations from $BV(\Omega)$}. 
Math. Mech. Complex Syst., \textbf{3} (2015), 83-100.

\bibitem{S} M. \v{S}ilhav\'y: \emph{The generals form of the relaxation of a purely interfacial energy for structured deformations}. Math. Mech. Complex Syst., \emph{to appear}.

\end{thebibliography}
\end{document}